\documentclass{amsart}%
\usepackage{amsmath}
\usepackage{amsfonts}
\usepackage{amssymb}
\usepackage{graphicx}%
\providecommand{\U}[1]{\protect\rule{.1in}{.1in}}
\newtheorem{theorem}{Theorem}
\theoremstyle{plain}

\newtheorem{corollary}{Corollary}

\newtheorem{definition}{Definition}

\newtheorem{lemma}{Lemma}

\newtheorem{proposition}{Proposition}
\newtheorem{remark}{Remark}

\numberwithin{equation}{section}
\begin{document}
\title{A Note On Inner Quasidiagonal C*-Algebras }
\author{Qihui Li}
\address{Department of Mathematics, East China University of Science and Technology,
Meilong Road 130, 200237 Shanghai, P.R. China.}
\email{qihui\_li@126.com}
\thanks{The research of the first author is partially supported by National Natural
Science Foundation of China.}

\begin{abstract}
In the paper, we give two new characterizations of separable inner
quasidiagonal C*-algebras. Base on these characterizations, we show that a
unital full free product of two inner quasidiagonal C*-algebras is inner
quasidiagonal again. As an application, we show that a unital full free
product of two inner quasidiagoanl C*-algebras with amalgmation over a full
matrix algebra is inner quasidiagonal. Meanwhile, we conclude that a unital
full free product of two AF algebras with amalgamation over a
finite-dimensional C*-algebra is inner quasidiagonal if there are faithful
tracial states on each of these two AF algebras such that the restrictions on
the common subalgebra agree.

\end{abstract}
\maketitle

\section{Introduction}

Quasidiagonal (QD) C*-algebras have now been studied for more than 30 years.
Voiculescu \cite{V} give a characterization of quasidiagonal C*-algebras as following:

\begin{definition}
A C*-algebra $\mathcal{A}$ is quasidiagonal if, for every $x_{1},\cdots
,x_{n}\in\mathcal{A}$ and $\varepsilon>0,$ there is a representation $\pi$ of
$\mathcal{A}$ on a Hilbert space $\mathcal{H}$, and a finite-rank projection
$P\in\mathcal{B}\left(  \mathcal{H}\right)  $ such that $\left\Vert
P\pi\left(  x_{i}\right)  -\pi\left(  x_{i}\right)  P\right\Vert
<\varepsilon,\left\Vert P\pi\left(  x_{i}\right)  P\right\Vert >\left\Vert
x_{i}\right\Vert -\varepsilon$ for $1\leq i\leq n.$
\end{definition}

Voiculescu showed that $\mathcal{A}$ is QD if and only if $\pi\left(
\mathcal{A}\right)  $ is a quasidiagonal set of operators for a faithful
essential representation $\pi$ of $\mathcal{A}$. In \cite{BK1}, we know that
all separable QD C*-algebras are Blackadar and Kirchberg's MF algebras. It is
well known that the reduced free group C*-algebra $C_{r}^{\ast}\left(
F_{2}\right)  $ is not QD. Haagerup and Thorbj$\phi$rnsen showed that
$C_{r}^{\ast}\left(  F_{2}\right)  $ is MF \cite{HT}. This implies that the
family of all separable QD C*-algebras are strictly contained in the set of MF C*-algebras.

The concept of MF algebras was first introduced by Blackadar and Kirchberg in
\cite{BK1}. Many properties of MF algebras were discussed in \cite{BK1}. In
the same article, Blackadar and Kirchberg study NF algebras and strong NF
algebras as well. A separable C*-algebra is a strong NF algebra if it can be
written as a generalized inductive limit of a sequential inductive system of
finite-dimensional C*-algebras in which the connecting maps are complete order
embedding and asymptotically multiplicative in the sense of \cite{BK1}. An NF
algebra is a C*-algebra which can be written as the generalized inductive
limit of such system, where the connecting maps are only required to be
completely positive contractions. It was shown that a separable C*-algebra is
an NF algebra if and only if it is nuclear and quasidiagonal. Whether the
class of NF algebra is distinct from the class of strong NF algebras? For
solving this question, Blackadar and Kirchberg introduce the concept of inner
quasidiagonal by slightly modifying Voiculescu's characterization of
quasidiagonal C*-algebras:

\begin{definition}
(\cite{BK2}) A C*-algebra $\mathcal{A}$ is inner quasidiagonal if, for every
$x_{1},\cdots,x_{n}\in\mathcal{A}$ and $\varepsilon>0,$ there is a
representation $\pi$ of $\mathcal{A}$ on a Hilbert space $\mathcal{H}$, and a
finite-rank projection $P\in\mathcal{\pi}\left(  \mathcal{A}\right)
^{\prime\prime}$ such that $\left\Vert P\pi\left(  x_{i}\right)  -\pi\left(
x_{i}\right)  P\right\Vert <\varepsilon,\left\Vert P\pi\left(  x_{i}\right)
P\right\Vert >\left\Vert x_{i}\right\Vert -\varepsilon$ for $1\leq i\leq n.$
\end{definition}

It was shown that a separable C*-algebra is a strong NF algebra if and only if
it is nuclear and inner quasidiagonal \cite{BK2}. Blackadar and Kirchberg also
gave examples of separable nuclear C*-algebras which are quasidiagonal but not
inner quasidiagonal, hence of NF algebras which are not strong NF. Therefore,
the preceding question has been solved.

In this note, we are interested in the question of whether the unital full
free products of inner QD C*-algebras are inner QD again. Note that every RFD
C*-algebra is inner QD \cite{BK2}. We have known that a unital full free
products of two RFD C*-algebras is RFD \cite{RT}. Similar result holds for
unital QD C*-algebras \cite{B}. Based on these results and the relationship
among RFD C*-algebras, inner QD C*-algebras and QD C*-algebras, it is natural
to ask whether the same things will happen when we consider inner QD
C*-algebras. In this note we will show that a unital full free product of two
unital inner QD C*-algebra is inner again. As an application, we will consider
the unital full free products of two inner QD C*-algebras with amalgamation
over finite-dimensional C*-algebras.

All C*-algebras in this note are unital and separable. A brief overview of
this paper is as follows. In Section 2, we fix some notation and give two new
characterizations of inner QD C*-algebras. Section 3 is devoted to results on
the unital full free products of two unital inner QD C*-algebras. We first
consider unital full free products of unital inner QD C*-algebras. As an
application, we show that a unital full free product of two inner
quasidiagonal C*-algebras with amalgamation over a full matrix algebra is
inner quasidiagonal. Meanwhile, we conclude that a unital full free product of
two AF algebras with amalgamation over a finite-dimensional C*-algebra is
inner quasidiagonal if there are faithful tracial states on each of these two
AF algebras such that the restrictions on the common subalgebra agree.

\section{Inner Quasidiagonal C*-algebras}

We denote the set of all bounded operators on $\mathcal{H}$ by $\mathcal{B}%
\left(  \mathcal{H}\right)  $.

Suppose $\{\mathcal{M}_{k_{n}}(\mathbb{C})\}_{n=1}^{\infty}$ is a sequence of
complex matrix algebras. We introduce the C*-direct product $\prod
_{m=1}^{\infty}\mathcal{M}_{k_{m}}(\mathbb{C)}$ of $\{\mathcal{M}_{k_{n}%
}(\mathbb{C})\}_{n=1}^{\infty}$ as follows:
\[
\prod_{n=1}^{\infty}\mathcal{M}_{k_{n}}(\mathbb{C})=\{(Y_{n})_{n=1}^{\infty
}\ |\ \forall\ n\geq1,\ Y_{n}\in\mathcal{M}_{k_{n}}(C)\ \text{ and
}\ \left\Vert (Y_{n})_{n=1}^{\infty}\right\Vert =\sup_{n\geq1}\Vert Y_{n}%
\Vert<\infty\}.
\]
Furthermore, we can introduce a norm-closed two sided ideal in $\prod
_{n=1}^{\infty}\mathcal{M}_{k_{n}}(\mathbb{C})$ as follows:
\[
\overset{\infty}{\underset{n=1}{\sum}}\mathcal{M}_{k_{n}}(\mathbb{C})=\left\{
\left(  Y_{n}\right)  _{n=1}^{\infty}\in\prod_{n=1}^{\infty}\mathcal{M}%
_{k_{n}}(\mathbb{C}):\lim\limits_{n\rightarrow\infty}\left\Vert Y_{n}%
\right\Vert =0\right\}  .
\]
Let $\pi$ be the quotient map from $\prod_{n=1}^{\infty}\mathcal{M}_{k_{n}%
}(\mathbb{C})$ to $\prod_{n=1}^{\infty}\mathcal{M}_{k_{n}}(\mathbb{C}%
)/\overset{\infty}{\underset{n=1}{\sum}}\mathcal{M}_{k_{n}}(\mathbb{C})$.
Then
\[
\prod_{n=1}^{\infty}\mathcal{M}_{k_{n}}(\mathbb{C})/\overset{\infty}%
{\underset{n=1}{\sum}}\mathcal{M}_{k_{n}}(\mathbb{C})
\]
is a unital C*-algebra. If we denote $\pi\left(  \left(  Y_{n}\right)
_{n=1}^{\infty}\right)  $ by $\left[  \left(  Y_{n}\right)  _{n}\right]  $,
then
\[
\left\Vert \left[  \left(  Y_{n}\right)  _{n}\right]  \right\Vert
=\underset{n\rightarrow\mathcal{1}}{\lim\sup}\left\Vert Y_{n}\right\Vert
\leq\sup_{n}\left\Vert Y_{n}\right\Vert =\left\Vert \left(  Y_{n}\right)
_{n}\right\Vert \in\prod_{n=1}^{\infty}\mathcal{M}_{k_{n}}(\mathbb{C})
\]

Recall that a C*-algebra is residually finite-dimensional (RFD) if it has a
separating family of finite-dimensional representations. If a separable C*-
algebra $\mathcal{A}$ can be embedded into $%
{\textstyle\prod\limits_{k}}
\mathcal{M}_{n_{k}}\left(  \mathbb{C}\right)  /\sum_{k}\mathcal{M}_{n_{k}%
}\left(  \mathbb{C}\right)  $ for a sequence of positive integers $\left\{
n_{k}\right\}  _{k=1}^{\mathcal{1}},$ then $\mathcal{A}$ is called an MF
algebra. Many properties of MF algebras were discussed in \cite{BK1}. Note
that the family of all RFD C*-algebras is strictly contained in the family of
all inner QD C*-algebras, and all QD C*-algebras are MF C*-algebras.

Continuing the study of generalized inductive limits of finite-dimensional
C*-algebras, Blackadar and Kirchberg define a refined notion of
quasidiagonality for C*-algebras, called inner quasidiagonality. A cleaner
alternative definition of inner quasidiagonality can be given using the socle
of the bidual.

\begin{definition}
If $\mathcal{B}$ is a C*-algebra, then a projection $p\in\mathcal{B}$ is in
the socle if $p\mathcal{B}p$ is finite-dimensional. Denote the set of the
socle in $\mathcal{B}$ by $socle\left(  \mathcal{B}\right)  $
\end{definition}

\begin{theorem}
\label{3}(\cite{BO}) A separable C*-algebra $\mathcal{A}$ is inner QD if there
are projections $p_{n}\in\mathcal{A}^{\ast\ast}$ such that

\begin{enumerate}
\item $\left\Vert \left[  p_{n},a\right]  \right\Vert \longrightarrow0$ for
all $a\in\mathcal{A\subseteq A}^{\ast\ast},$

\item $\left\Vert a\right\Vert =\lim\left\Vert p_{n}ap_{n}\right\Vert $ for
all $a\in\mathcal{A}$ and

\item $p_{n}\in socle\left(  \mathcal{A}^{\ast\ast}\right)  $ for every $n.$
\end{enumerate}
\end{theorem}

\begin{theorem}
\label{11a}(\cite{BK2}, Proposition 3.7.) Let $\mathcal{A}$ be a separable
C*-algebra. Then $\mathcal{A}$ is inner QD if and only if there is a sequence
of irreducible representation $\left\{  \pi_{n}\right\}  $ of $\mathcal{A}$ on
Hilbert space $\mathcal{H}_{n},$and finite-rank projection $p_{n}%
\in\mathcal{B}\left(  \mathcal{H}_{n}\right)  $, such that $\left\Vert \left[
p_{n},\pi_{n}\left(  x\right)  \right]  \right\Vert \longrightarrow0$ and
$\lim\sup\left\Vert p_{n}\pi_{n}\left(  x\right)  p_{n}\right\Vert =\left\Vert
x\right\Vert $ for all $x\in\mathcal{A}$.
\end{theorem}

The principal shortcoming of the definition of inner QD C*-algebra is that it
is often difficult to determine directly whether a C*-algebra is inner QD, the
following result for separable case is much easier for checking.

\begin{theorem}
\label{11b}(\cite{BK3}) A separable C*-algebra is inner QD if and only if it
has a separating family of quasidiagonal irreducible representations.
\end{theorem}

Let $\pi:\mathcal{B}\left(  \mathcal{H}\right)  \rightarrow\mathcal{Q}\left(
\mathcal{H}\right)  $ be the canonical mapping onto the Calkin algebra and
$\mathcal{A}$ is a unital C*-algebra. Suppose $\varphi:\mathcal{A\rightarrow
B}\left(  \mathcal{H}\right)  $ is a unital completely positive map then we
say that $\varphi$ is a representation modulo the compacts if $\pi\circ
\varphi:\mathcal{A\rightarrow Q}\left(  \mathcal{H}\right)  $ is a
*-homomorphism. If $\pi\circ\varphi$ is injective then we say that $\varphi$
is a faithful representation modulo the compacts.

For an MF C*-algebra, we are able to embed it into $%
{\textstyle\prod\limits_{k}}
\mathcal{M}_{n_{k}}\left(  \mathbb{C}\right)  /\sum_{k}\mathcal{M}_{n_{k}%
}\left(  \mathbb{C}\right)  $ for a sequence of positive integers $\left\{
n_{k}\right\}  _{k=1}^{\mathcal{1}}.$ For an RFD C*-algebra, we can embed it
into $%
{\textstyle\prod\limits_{k}}
\mathcal{M}_{n_{k}}\left(  \mathbb{C}\right)  $. Meanwhile, for a QD
C*-algebra, we can not only embed it into $%
{\textstyle\prod\limits_{k}}
\mathcal{M}_{n_{k}}\left(  \mathbb{C}\right)  /\sum_{k}\mathcal{M}_{n_{k}%
}\left(  \mathbb{C}\right)  ,$ but also lift this embedding to a faithful
representation into $%
{\displaystyle\prod}
\mathcal{M}_{k_{m}}\left(  \mathbb{C}\right)  $ modulo the compacts. Whether
there is a similar characterization for the inner QD C*-algebras? We will
answer this question in the following theorem.

The following lemma is a well-known result about completely positive map. We
use c.p. to abbreviate "completely positive", u.c.p. for "unital completely
positive" and c.c.p. for "contractive completely positive".

\begin{lemma}
\label{11}(Stinespring) Let $\mathcal{A}$ be a unital C*-algebra and
$\varphi:\mathcal{A\longrightarrow B}\left(  \mathcal{H}\right)  $ be a c.p..
map. Then, there exist a Hilbert space $\widehat{\mathcal{H}}$, a
*-representation $\pi_{\varphi}:\mathcal{A\longrightarrow B}\left(
\mathcal{H}\right)  $ and an operator $V:\mathcal{H\longrightarrow}%
\widehat{\mathcal{H}}$ such that
\[
\varphi\left(  a\right)  =V^{\ast}\pi_{\varphi}\left(  a\right)  V
\]
for every $a\in\mathcal{A}$. In particular, $\left\Vert \varphi\right\Vert
=\left\Vert V^{\ast}V\right\Vert =\left\Vert \varphi\left(  1\right)
\right\Vert .$
\end{lemma}

We call the triplet $\left(  \pi_{\varphi},\widehat{\mathcal{H}},V\right)  $
in preceding lemma a Stinespring dilation of $\varphi.$ When $\varphi$ is
unital, $V^{\ast}V=\varphi\left(  I\right)  =I,$ and hence $V$ is an isometry.
So in this case we may assume that $V$ is a projection $P$ and $\varphi\left(
a\right)  =P\pi_{\varphi}\left(  a\right)  |_{\mathcal{H}}$. In general there
could be many different Stinespring dilations, but we may always assume that a
dilation $\left(  \pi_{\varphi},\widehat{\mathcal{H}},V\right)  $ is minimal
in the sense that $\pi_{\varphi}\left(  \mathcal{A}\right)  V\mathcal{H}$ is
dense in $\widehat{\mathcal{H}}.$ We know that, under this minimality
condition, a Stinespring dilation is unique up to unitary equivalence. Note
that if $\left(  \pi_{\varphi},\widehat{\mathcal{H}},V\right)  $ is minimal
Stinespring dilation of $\varphi:\mathcal{A\longrightarrow B}\left(
\mathcal{H}\right)  $, then there exists a *-homomorphism $\rho:\varphi\left(
\mathcal{A}\right)  ^{\prime}\longrightarrow\pi_{\varphi}\left(
\mathcal{A}\right)  ^{\prime}\subseteq\mathcal{B}\left(  \widehat{\mathcal{H}%
}\right)  $ such that $\varphi\left(  a\right)  x=V^{\ast}\pi_{\varphi}\left(
a\right)  \rho\left(  x\right)  V$ for every $a\in\mathcal{A}$ and
$x\in\varphi\left(  \mathcal{A}\right)  ^{\prime},$ it implies that the
commutant $\varphi\left(  \mathcal{A}\right)  ^{\prime}\subseteq
\mathcal{B}\left(  \mathcal{H}\right)  $ also lifts to $\mathcal{B}\left(
\widehat{\mathcal{H}}\right)  .$

\begin{lemma}
\label{11c}Let $\mathcal{A}$ be a unital C*-algebra and $\varphi
:\mathcal{A\longrightarrow M}_{n}\left(  \mathbb{C}\right)  $ be a surjective
u.c.p. map. Suppose $\left(  \pi_{\varphi},\widehat{\mathcal{H}},P\right)  $
is a minimal Stinespring dilation of $\varphi$ where $P$ is a projection in
$\mathcal{B}\left(  \widehat{\mathcal{H}}\right)  .$ Then the *-homomorphism
$\rho:\varphi\left(  \mathcal{A}\right)  ^{\prime}\longrightarrow\pi_{\varphi
}\left(  \mathcal{A}\right)  ^{\prime}\subseteq\mathcal{B}\left(
\widehat{\mathcal{H}}\right)  $ is unital, i.e.%
\[
\rho\left(  \mathbb{C}I\right)  =\mathbb{C}I\subseteq\pi_{\varphi}\left(
\mathcal{A}\right)  ^{\prime}%
\]

\end{lemma}

\begin{proof}
Since $\varphi$ is surjective, $\varphi\left(  \mathcal{A}\right)  ^{\prime
}=\mathbb{C}I$. Note $\rho$ is a *-homomorphism, it is easy to check that
\[
\left(  I-P\right)  \rho\left(  \alpha I\right)  P=P\rho\left(  \alpha
I\right)  \left(  I-P\right)  =0.
\]
and $\left(  I-P\right)  \rho\left(  I\right)  \left(  I-P\right)  $ is a
projection. We know $\left(  \pi_{\varphi},\widehat{\mathcal{H}},P\right)  $
is a minimal Stinespring dilation, then $\pi_{\varphi}\left(  \mathcal{A}%
\right)  ^{\prime}$ has no proper projection bigger than $P.$ It implies
$\left(  I-P\right)  \rho\left(  I\right)  \left(  I-P\right)  =0,$ i.e.
$\rho\left(  I\right)  =I$. Hence $\rho\left(  \mathbb{C}I\right)
=\mathbb{C}I\subseteq\pi_{\varphi}\left(  \mathcal{A}\right)  ^{\prime}.$
\end{proof}

Now, we are ready to give a new characterization of inner QD C*-algebras.

\begin{theorem}
Suppose $\mathcal{A}$ is a unital C*-algebra. Then $\mathcal{A}$ is inner QD
if and only if there is a faithful representation modulo compacts
$\Phi:\mathcal{A\longrightarrow}$ $\Pi\mathcal{M}_{k_{n}}\left(
\mathbb{C}\right)  $ for a sequence $\left\{  k_{n}\right\}  $ of integers
such that the u.c.p. maps $\varphi_{n}:\mathcal{A\longrightarrow M}_{k_{n}%
}\left(  \mathbb{C}\right)  $ is surjective for every $n$ and the
*-homomorphism%
\[
\rho:\varphi_{n}\left(  \mathcal{A}\right)  ^{\prime}\longrightarrow
\pi_{\varphi_{n}}\left(  \mathcal{A}\right)  ^{\prime}\subseteq\mathcal{B}%
\left(  \widehat{\mathcal{H}}_{n}\right)
\]
is surjective where $\left(  \pi_{\varphi_{n}},\widehat{\mathcal{H}}_{n}%
,p_{n}\right)  $ is a minimal Stinespring dilation of $\varphi_{n}.$
\end{theorem}

\begin{proof}
$\left(  \Longrightarrow\right)  $ Suppose $\mathcal{A}$ is inner QD. Then, by
applying Theorem \ref{11a}, we can find sequences of irreducible
representations $\left\{  \pi_{n}\right\}  $ and finite projections $\left\{
p_{n}\right\}  $ where $p_{n}\in\pi_{n}\left(  \mathcal{A}\right)
^{\prime\prime}$ such that $\Phi:\mathcal{A\longrightarrow}$ $\Pi p_{n}\pi
_{n}\left(  \mathcal{A}\right)  p_{n}$ is a faithful representation modulo
compacts. Meanwhile, we have $p_{n}\pi_{n}\left(  \mathcal{A}\right)
p_{n}\cong\mathcal{M}_{k_{n}}\left(  \mathbb{C}\right)  $ for some integer
$k_{n}$ and $\left(  \pi_{n}\left(  \mathcal{A}\right)  \right)  ^{\prime
}=\mathbb{C}I$ since $\pi_{n}$ is irreducible. Define
\[
\varphi_{n}:\mathcal{A\longrightarrow}p_{n}\pi_{n}\left(  \mathcal{A}\right)
p_{n}\cong\mathcal{M}_{k_{n}}\left(  \mathbb{C}\right)
\]
by $\varphi_{n}\left(  a\right)  =p_{n}\pi_{n}\left(  a\right)  p_{n}$. Then
$\varphi_{n}$ is u.c.p. and surjective for every $n.$ Note $\left(  \pi
_{n},\mathcal{H}_{n},p_{n}\right)  $ is a minimal Stinespring dilation of
$\varphi_{n}$ since $\pi_{n}$ is irreducible. Therefore the *-homomorphism
\[
\rho:\varphi_{n}\left(  \mathcal{A}\right)  ^{\prime}\longrightarrow\pi
_{n}\left(  \mathcal{A}\right)  ^{\prime}\subseteq\mathcal{B}\left(
\mathcal{H}_{n}\right)
\]
is surjective by Lemma \ref{11c} and the fact that $\pi\left(  \mathcal{A}%
\right)  ^{\prime}=\mathbb{C}I.$

$\left(  \Longleftarrow\right)  $ Suppose there is a faithful representation
modulo compacts $\Phi:\mathcal{A\longrightarrow}$ $\Pi\mathcal{M}_{k_{n}%
}\left(  \mathbb{C}\right)  $ for a sequence $\left\{  k_{n}\right\}  $ such
that the u.c.p. maps $\varphi_{n}:\mathcal{A\longrightarrow M}_{k_{n}}\left(
\mathbb{C}\right)  $ is surjective and the *-homomorphism
\[
\rho:\varphi_{n}\left(  \mathcal{A}\right)  ^{\prime}\longrightarrow
\pi_{\varphi_{n}}\left(  \mathcal{A}\right)  ^{\prime}\subseteq\mathcal{B}%
\left(  \widehat{\mathcal{H}}_{n}\right)
\]
is surjective where $\left(  \pi_{\varphi_{n}},\widehat{\mathcal{H}}_{n}%
,p_{n}\right)  $ is a minimal Stinespring dilation of $\varphi_{n}.$ Then
\[
\rho\left(  \varphi_{n}\left(  \mathcal{A}\right)  ^{\prime}\right)
=\rho\left(  \mathbb{C}I\right)  =\mathbb{C}I
\]
by Lemma \ref{11c}. It implies that $\pi_{\varphi_{n}}\left(  \mathcal{A}%
\right)  ^{\prime}=\mathbb{C}I$ since $\rho$ is surjective. Hence
$\pi_{\varphi_{n}}$ is irreducible and $p_{n}\in\pi_{\varphi_{n}}\left(
\mathcal{A}\right)  ^{\prime\prime}.$ So, for these irreducible representation
$\left\{  \pi_{\varphi_{n}}\right\}  $ of $\mathcal{A}$ on Hilbert space
$\widehat{\mathcal{H}}_{n}$ and finite-rank projection $p_{n}\in\pi
_{\varphi_{n}}\left(  \mathcal{A}\right)  ^{\prime\prime}$, we have
$\left\Vert \left[  p_{n},\pi_{\varphi_{n}}\left(  x\right)  \right]
\right\Vert \longrightarrow0$ and $\lim\sup\left\Vert p_{n}\pi_{\varphi_{n}%
}\left(  x\right)  p_{n}\right\Vert =\left\Vert x\right\Vert $ for all
$x\in\mathcal{A}$. It implies that $\mathcal{A}$ is inner QD by Theorem
\ref{11a}.
\end{proof}

Suppose $\mathcal{A}$ is a unital C*-algebra and $p\in socle\left(
\mathcal{A}^{\ast\ast}\right)  .$ Define
\[
\mathcal{A}_{p}=\left\{  a\in\mathcal{A}:\left[  a,p\right]  =0\right\}  .
\]
Then we have the following few lemmas.

\begin{lemma}
\label{12a}(\cite{BK2}, Corollary 3.5.) Let $p\in socle\left(  \mathcal{A}%
^{\ast\ast}\right)  .$ Then $d\left(  a,\mathcal{A}_{p}\right)  =\left\Vert
\left[  a,p\right]  \right\Vert $for all $a\in\mathcal{A}$.
\end{lemma}

\begin{lemma}
\label{12b}(\cite{BK2}, Proposition 3.4.) Let $\mathcal{A}$ be a C*-algebra,
and $p\in socle\left(  \mathcal{A}^{\ast\ast}\right)  .$ Then

\begin{enumerate}
\item $p\mathcal{A}_{p}=p\mathcal{A}_{p}p=p\mathcal{A}^{\ast\ast
}p=p\mathcal{A}p$

\item The weak closure of $\mathcal{A}_{p}$ in $\mathcal{A}^{\ast\ast}$ is
$p\mathcal{A}^{\ast\ast}p+\left(  1-p\right)  \mathcal{A}^{\ast\ast}(1-p).$
\end{enumerate}
\end{lemma}

\begin{lemma}
\label{12ee}Let $\mathcal{A}$ be a C*-algebra, $p_{1},\cdots,p_{k}\in
socle\left(  \mathcal{A}^{\ast\ast}\right)  $ with $p_{1}\leq\cdots\leq
p_{k}.$Then
\begin{align*}
p_{k}\left(  \cap_{i=0}^{k}\mathcal{A}_{p_{i}}\right)   &  =p_{k}\left(
\cap_{i=0}^{k}\mathcal{A}_{p_{i}}^{\ast\ast}\right) \\
&  =p_{0}\mathcal{A}p_{0}\oplus\left(  p_{1}-p_{0}\right)  \mathcal{A}\left(
p_{1}-p_{0}\right)  \oplus\cdots\oplus\left(  p_{k}-p_{k-1}\right)
\mathcal{A}\left(  p_{k}-p_{k-1}\right)
\end{align*}

\end{lemma}

\begin{proof}
We only prove the case when $k=2.$ Since $p_{1}\leq p_{2}\in socle\left(
\mathcal{A}^{\ast\ast}\right)  ,$ we have
\[
p_{2}\left(  \mathcal{A}_{p_{1}}^{\ast\ast}\cap\mathcal{A}_{p_{2}}^{\ast\ast
}\right)  =p_{1}\mathcal{A}p_{1}+\left(  p_{2}-p_{1}\right)  \mathcal{A}%
\left(  p_{2}-p_{1}\right)
\]
by Lemma $\ref{12b}$. So it is obvious that
\[
p_{2}\left(  \mathcal{A}_{p_{1}}\cap\mathcal{A}_{p_{2}}\right)  \subseteq
p_{1}\mathcal{A}p_{1}+\left(  p_{2}-p_{1}\right)  \mathcal{A}\left(
p_{2}-p_{1}\right)  .
\]
Meanwhile,
\begin{align*}
p_{2}\mathcal{A}_{p_{2}}  &  =p_{2}\mathcal{A}p_{2}=p_{2}\mathcal{A}^{\ast
\ast}p_{2}\\
&  \supseteq p_{1}\mathcal{A}p_{1}+\left(  p_{2}-p_{1}\right)  \mathcal{A}%
\left(  p_{2}-p_{1}\right)
\end{align*}
by Lemma $\ref{12b}$ and the fact that $p_{1}\leq p_{2}\in socle\left(
\mathcal{A}^{\ast\ast}\right)  .$ Then for every $b_{1}\in p_{1}%
\mathcal{A}p_{1}$ and $b_{2}\in\left(  p_{2}-p_{1}\right)  \mathcal{A}\left(
p_{2}-p_{1}\right)  ,$there is $a\in\mathcal{A}_{p_{2}}$ such that
$p_{2}a=b_{1}+b_{2}.$ Hence $a\in\mathcal{A}_{p_{1}}\cap\mathcal{A}_{p_{2}}$.
It implies that%
\[
p_{2}\left(  \mathcal{A}_{p_{1}}\cap\mathcal{A}_{p_{2}}\right)  \supseteq
p_{1}\mathcal{A}p_{1}+\left(  p_{2}-p_{1}\right)  \mathcal{A}\left(
p_{2}-p_{1}\right)  .
\]

This completes the proof.
\end{proof}

\begin{lemma}
\label{12ef}Let $\mathcal{A}$ be a C*-algebra, $\left\{  p_{n}\right\}  $ be a
sequence of projection in $socle\left(  \mathcal{A}^{\ast\ast}\right)  $ with
$p_{0}\leq p_{1}\leq\cdots$ and $p_{i}\overset{s.o.t}{\longrightarrow}I$
(strong operator topology). Then%
\[
\cap_{i=0}^{\mathcal{1}}\mathcal{A}_{p_{i}}=p_{0}\mathcal{A}p_{0}\oplus
\oplus_{i=1}^{\mathcal{1}}\left[  \left(  p_{i}-p_{i-1}\right)  \mathcal{A}%
\left(  p_{i}-p_{i-1}\right)  \right]  \subseteq\mathcal{A}%
\]

\end{lemma}

\begin{proof}
By Lemma \ref{12ee} and the fact that $p_{1}\leq p_{2}\leq\cdots$ with
$p_{i}\in socle\left(  \mathcal{A}^{\ast\ast}\right)  ,$ we have
\begin{align*}
p_{k}\left(  \cap_{i=0}^{\mathcal{1}}\mathcal{A}_{p_{i}}\right)   &
=p_{k}\left(  \cap_{i=0}^{k}\mathcal{A}_{p_{i}}\right)  =p_{k}\left(
\cap_{i=0}^{k}\mathcal{A}_{p_{i}}^{\ast\ast}\right) \\
&  =p_{0}\mathcal{A}p_{0}\oplus\left(  p_{1}-p_{0}\right)  \mathcal{A}\left(
p_{1}-p_{0}\right)  \oplus\cdots\oplus\left(  p_{k}-p_{k-1}\right)
\mathcal{A}\left(  p_{k}-p_{k-1}\right)
\end{align*}
for every $k$. Therefore
\begin{align*}
\cap_{i=0}^{\mathcal{1}}\mathcal{A}_{p_{i}}  &  =p_{0}\left(  \cap
_{i=0}^{\mathcal{1}}\mathcal{A}_{p_{i}}\right)  \oplus\oplus_{i=1}%
^{\mathcal{1}}\left(  p_{i}-p_{i-1}\right)  \left(  \cap_{k=0}^{\mathcal{1}%
}\mathcal{A}_{p_{k}}\right) \\
&  =p_{0}\left(  \cap_{i=0}^{\mathcal{1}}\mathcal{A}_{p_{i}}\right)
\oplus\oplus_{i=1}^{\mathcal{1}}\left(  p_{i}-p_{i-1}\right)  \left(
\cap_{k=0}^{i}\mathcal{A}_{p_{k}}\right) \\
&  =p_{0}\mathcal{A}p_{0}\oplus\oplus_{i=1}^{\mathcal{1}}\left[  \left(
p_{i}-p_{i-1}\right)  \mathcal{A}\left(  p_{i}-p_{i-1}\right)  \right]
\end{align*}

\end{proof}

\begin{remark}
\label{12e}Suppose $\mathcal{A}$ is a unital inner QD C*-algebras, then there
is sequence $\left\{  p_{n}\right\}  $ of projections in $socle\left(
\mathcal{A}^{\ast\ast}\right)  $ such that $\left\Vert \left[  p_{n},a\right]
\right\Vert \longrightarrow0$ for all $a\in\mathcal{A\subseteq A}^{\ast\ast}$
and $\left\Vert a\right\Vert =\lim\left\Vert p_{n}ap_{n}\right\Vert $ for all
$a\in\mathcal{A}$ by Theorem \ref{3}. Therefore we can define a sequence of
u.c.p maps $\varphi_{n}:\mathcal{A\longrightarrow}p_{n}\mathcal{A}^{\ast\ast
}p_{n}$ by compression. It is obvious that $\mathcal{A}_{p_{n}}=\mathcal{M}%
_{\varphi_{n}}$ where $\mathcal{M}_{\varphi_{n}}$ is the multiplicative domain
of $\varphi_{n}$ and $\left\Vert a\right\Vert =\lim\left\Vert \varphi
_{n}\left(  a\right)  \right\Vert $, $d\left(  a,\mathcal{M}_{\varphi_{n}%
}\right)  \longrightarrow0$ for all $a\in\mathcal{A}$ by Lemma \ref{12a}$.$
Actually, this is a sufficient condition for a given C*-algebra to be an inner
QD C*-algebra.
\end{remark}

\begin{theorem}
\label{12d}(\cite{BO}) $\mathcal{A}$ is inner QD if and only if there is a
sequence of c.c.p. maps $\varphi_{n}:\mathcal{A\longrightarrow M}_{k_{n}%
}\left(  \mathbb{C}\right)  $ such that $\left\Vert a\right\Vert
=\lim\left\Vert \varphi_{n}\left(  a\right)  \right\Vert $ and $d\left(
a,\mathcal{M}_{\varphi_{n}}\right)  \longrightarrow0$ for all $a\in
\mathcal{A}$, where $\mathcal{M}_{\varphi_{n}}$ is the multiplicative domain
of $\varphi_{n}.$
\end{theorem}

Now, we are ready to give another characterization of unital inner QD C*-algebras.

\begin{theorem}
\label{ttt}Suppose $\mathcal{A}$ is a unital separable C*-algebra. Then
$\mathcal{A}$ is inner QD if and only if there is a sequence of unital RFD
C*-subalgebra $\left\{  \mathcal{A}_{n}\right\}  _{n=1}^{\mathcal{1}}$ of
$\mathcal{A}$ such that $\cup_{n=1}^{\mathcal{1}}\mathcal{A}_{n}$ is norm
dense in $\mathcal{A}$.

\begin{proof}
$\left(  \Longrightarrow\right)  $ Suppose $\mathcal{F}\subseteq\mathcal{A}$
is a finite subset and $\varepsilon>0$. Let
\[
\left\{  1\right\}  \cup\mathcal{F}\subseteq\mathcal{F}_{1}\subseteq
\mathcal{F}_{2}\mathcal{\subseteq\cdots}%
\]
be the sequence of finite subsets of $\mathcal{A}$ such that $\overline
{\cup_{i}\mathcal{F}_{i}}=\mathcal{A}.$ Then$,$ from Remark \ref{12e} and
Theorem \ref{12d}, we can find
\begin{align*}
P_{0} &  \leq P_{1}\leq P_{2}\leq\cdots\text{ \ }\\
\text{with }P_{i} &  \in socle\left(  \mathcal{A}^{\ast\ast}\right)  \text{
and }P_{i}\overset{s.o.t}{\longrightarrow}P\in\mathcal{A}^{\ast\ast}\text{
}(i\longrightarrow\mathcal{1})
\end{align*}
such that
\begin{align*}
d(a,\mathcal{A}_{P_{i}}) &  =\left\Vert [a,P_{i}]\right\Vert \\
&  =\left\Vert \left(  1-P_{i}\right)  aP_{i}+P_{i}a\left(  1-P_{i}\right)
\right\Vert <\frac{\varepsilon}{2\cdot2^{i+1}}\text{ }%
\end{align*}
and $\left\Vert P_{i}aP_{i}\right\Vert >\left\Vert a\right\Vert -\frac
{\varepsilon}{2^{i+1}}$ for every $a\in\mathcal{F}_{i}$ $(i\in\mathbb{N)}$.
Since $P_{i}\overset{s.o.t}{\longrightarrow}P\in\mathcal{A}^{\ast\ast}$ $($as
$i\longrightarrow\mathcal{1})$ and $P\geq P_{i},$ we have
\[
\left\Vert PaP\right\Vert \geq\left\Vert P_{i}aP_{i}\right\Vert \geq\left\Vert
a\right\Vert -\frac{\varepsilon}{2^{i+1}}\text{ for }\forall a\in\cup
_{i}\mathcal{F}_{i}\text{ and }i.
\]
It implies that $\left\Vert PaP\right\Vert =\left\Vert a\right\Vert $ for
$\forall a\in\overline{\cup_{i}\mathcal{F}_{i}}=\mathcal{A}$, therefore $P=I.$
Now let
\[
\mathcal{A}_{\varepsilon}=\cap_{i=0}^{\mathcal{1}}\mathcal{A}_{P_{i}}%
=P_{0}\mathcal{A}P_{0}\oplus\left(  P_{1}-P_{0}\right)  \mathcal{A}\left(
P_{1}-P_{0}\right)  \oplus\cdots
\]
by Lemma \ref{12ef}. So, for any $a\in\mathcal{F},$ let
\[
x=P_{0}aP_{0}+\left(  P_{1}-P_{0}\right)  a\left(  P_{1}-P_{0}\right)
+\cdots\in\mathcal{A}_{\varepsilon},
\]
we have
\begin{align*}
d(a,\mathcal{A}_{\varepsilon}) &  \leq\left\Vert a-x\right\Vert \\
&  =||P_{0}a\left(  P_{1}-P_{0}\right)  +P_{1}a\left(  P_{2}-P_{1}\right)
+\cdots\\
&  +\left(  P_{1}-P_{0}\right)  aP_{0}+\left(  P_{2}-P_{1}\right)
aP_{1}+\cdots||\\
&  \leq\left\Vert P_{0}a\left(  1-P_{0}\right)  \right\Vert \left\Vert
P_{1}-P_{0}\right\Vert +\cdots+\left\Vert P_{1}-P_{0}\right\Vert \left\Vert
\left(  1-P_{0}\right)  aP_{0}\right\Vert \\
&  <\sum_{i=0}^{\mathcal{1}}\frac{\varepsilon}{2\cdot2^{i+1}}+\sum
_{i=0}^{\mathcal{1}}\frac{\varepsilon}{2\cdot2^{i+1}}=\varepsilon.
\end{align*}
Note that $\mathcal{A}_{\varepsilon}$ is an RFD C*-subalgebras of
$\mathcal{A}$, hence we can find a sequence of unital RFD C*-subalgebra
$\left\{  \mathcal{A}_{n}\right\}  _{n=1}^{\mathcal{1}}$ of $\mathcal{A}$ such
that $\cup_{n=1}^{\mathcal{1}}\mathcal{A}_{n}$ is norm dense in $\mathcal{A}$.

$\left(  \Longleftarrow\right)  $ Suppose $\left\{  \mathcal{A}_{n}\right\}
_{n=1}^{\mathcal{1}}$ is a sequence of unital RFD C*-subalgebra in
$\mathcal{A}$ such that $\overline{\cup_{n=1}^{\mathcal{1}}\mathcal{A}_{n}%
}^{\left\Vert \cdot\right\Vert }=\mathcal{A}$, $\mathcal{F\subseteq A}$ is a
finite subset and $\varepsilon>0.$ Then there is an RFD C*-subalgebra
$\mathcal{A}_{n}$ and $b_{a}\in\mathcal{A}_{n}$ such that $\left\Vert
a-b_{a}\right\Vert <\frac{\varepsilon}{3}$ for every $a\in\mathcal{F}$. It
follows that%
\[
\left\Vert a\right\Vert -\frac{\varepsilon}{3}\leq\left\Vert b_{a}\right\Vert
.
\]
Since $\mathcal{A}_{n}$ is RFD, we can find a projection $P$ such that
$\Phi_{P}:\mathcal{A}_{n}\longrightarrow P\mathcal{A}_{n}P\subseteq
\mathcal{M}_{t}\left(  \mathbb{C}\right)  $ is a *-homorphism for some
$t\in\mathbb{C}$ and $\left\Vert \Phi_{P}\left(  b_{a}\right)  \right\Vert
\geq\left\Vert b_{a}\right\Vert -\frac{\varepsilon}{3}.$ Extending $\Phi_{P}$
to a u.c.p. map $\widetilde{\Phi_{P}}:$ $\mathcal{A\longrightarrow M}%
_{t}\left(  \mathbb{C}\right)  $ with $\mathcal{A}_{n}\subseteq\mathcal{M}%
_{\widetilde{\Phi_{P}}}$ and
\[
\left\Vert \widetilde{\Phi_{P}}\left(  b_{a}\right)  \right\Vert
\geq\left\Vert b_{a}\right\Vert -\frac{\varepsilon}{3}%
\]
where $\mathcal{M}_{\widetilde{\Phi_{P}}}$ is the multiplicative domain of
$\widetilde{\Phi_{P}}.$ Then
\[
\left\Vert \widetilde{\Phi_{P}}\left(  b_{a}\right)  \right\Vert =\left\Vert
\widetilde{\Phi_{P}}\left(  b_{a}-a\right)  +\widetilde{\Phi_{P}}\left(
a\right)  \right\Vert \leq\frac{\varepsilon}{3}+\left\Vert \widetilde{\Phi
_{P}}\left(  a\right)  \right\Vert \text{ for every }a\in\mathcal{F}\text{.}%
\]
So from above inequalities, we have
\begin{align*}
\left\Vert \widetilde{\Phi_{P}}\left(  a\right)  \right\Vert  &
\geq\left\Vert \widetilde{\Phi_{P}}\left(  b_{a}\right)  \right\Vert
-\frac{\varepsilon}{3}\\
&  \geq\left\Vert b_{a}\right\Vert -\frac{2\varepsilon}{3}\geq\left\Vert
a\right\Vert -\varepsilon\text{ for every }a\in\mathcal{F}%
\end{align*}
By the preceding discussion, for a finite subset $\mathcal{F}$ and
$\varepsilon>0$, there is a u.c.p map $\widetilde{\Phi_{P}}:\mathcal{A}%
\longrightarrow\mathcal{M}_{t}\left(  \mathbb{C}\right)  $ for some
$t\in\mathbb{N}$ such that $d\left(  a,\mathcal{M}_{\widetilde{\Phi_{P}}%
}\right)  <\varepsilon$ and $\left\Vert \widetilde{\Phi_{P}}\left(  a\right)
\right\Vert \geq\left\Vert a\right\Vert -\varepsilon$ for every $a\in
\mathcal{F}.$ So by Theorem \ref{12d}, $\mathcal{A}$ is inner.
\end{proof}
\end{theorem}

\section{Unital Full Free Products of Two Inner QD Algebras.}

In this section we will consider the question of whether the unital full free
products of inner QD C*-algebras are inner QD again. First, we need a lemma
for showing the main result in this section.

\begin{lemma}
(Theorem 3.2, \cite{RT})\label{13.4} Suppose $\mathcal{A}_{1}$ and
$\mathcal{A}_{2}$ are unital C*-algebras. Then the unital full free product
$\mathcal{A=\mathcal{A}}_{1}\mathcal{\ast}_{\mathbb{C}}\mathcal{A}_{2}$ is RFD
if and only if $\mathcal{A}_{1}$ and $\mathcal{A}_{2}$ are both RFD.
\end{lemma}

Now, we are ready to give the main result of this section.

\begin{theorem}
\label{12g}If $\mathcal{A}_{1}$ and $\mathcal{A}_{2}$ are both unital inner
quasidiagonal C*-algebras. Then $\mathcal{A}_{1}\ast_{\mathbb{C}}%
\mathcal{A}_{2}$ is inner QD.
\end{theorem}

\begin{proof}
Suppose $\tau$ is a fixed state on $\mathcal{A}_{1}\ast_{\mathbb{C}%
}\mathcal{A}_{2}$ and $\mathcal{F}$ is a finite subset of $\mathcal{A}_{1}%
\ast_{\mathbb{C}}\mathcal{A}_{2}.$ Let
\[
\mathcal{A}_{j}^{0}=\left\{  a\in\mathcal{A}_{j}:\tau\left(  a\right)
=0\right\}  ,j=1,2.
\]
No loss of generality, we may assume that every $b\in\mathcal{F}$ can be
decomposed into a finite sum with respect to $\tau,$ that is%
\[
b=\alpha_{0}I+\sum_{i_{1}\neq i_{2}\neq\cdots\neq i_{n}}a_{i_{1}}a_{i_{2}%
}\cdots a_{i_{n}}\text{ \ \ }\alpha_{0}\in\mathbb{C}\text{, }a_{i_{j}}%
\in\mathcal{A}_{i_{j}}^{0},i_{1}\neq i_{2}\neq\cdots\neq i_{n}%
\]
where $\mathcal{A}_{i_{j}}^{0}=\mathcal{A}_{1}^{0}$ or $\mathcal{A}_{2}^{0}.$
Denote by $\mathcal{F}_{0}^{j}$, $j=1,2,$ the set of such elements of
$\mathcal{A}_{j}^{0}$ which appear in the decomposition of elements from
$\mathcal{F}$. Then we can find an RFD C*-subalgebra $\mathcal{A}%
_{\varepsilon}^{j}$ of $\mathcal{A}_{j}$ for $j=1,2$ such that
\begin{equation}
d(a,\mathcal{A}_{\varepsilon}^{j})<\varepsilon\text{ for }\forall
a\in\mathcal{F}_{0}^{j},j=1,2. \tag{1}\label{a}%
\end{equation}

Let $b$ be an element in $\mathcal{F}$. No loss of generality, we may assume
that $b$ can be decomposed into the form
\[
\alpha I+\sum_{i=1}^{l}a_{1,1}^{i}a_{2,1}^{i}a_{1,2}^{i}a_{2,2}^{i}a_{1,3}%
^{i}\cdots a_{1,n_{i}}^{i},
\]
where
\[
\left\{  a_{1,1}^{i},a_{1,2}^{i},\cdots,a_{1,n_{i}}^{i}\right\}
\subseteq\mathcal{F}_{0}^{1}%
\]
and
\[
\left\{  a_{2,1}^{i},a_{2,2}^{i},\cdots,a_{2,n}^{i}\right\}  \subseteq
\mathcal{F}_{0}^{2}.
\]
Then, for any $a_{j,k}^{i},$ we can find $\widetilde{a_{j,k}^{i}}%
\in\mathcal{A}_{\varepsilon}^{j}$ such that $\left\Vert a_{j,k}^{i}%
-\widetilde{a_{j,k}^{i}}\right\Vert <\varepsilon$ by (\ref{a}). Note that
\[
\left\vert \tau\left(  a_{j,k}^{i}-\widetilde{a_{j,k}^{i}}\right)  \right\vert
=\left\vert \tau\left(  a_{j,k}^{i}\right)  -\tau\left(  \widetilde
{a_{j,k}^{i}}\right)  \right\vert =\left\vert \tau\left(  \widetilde
{a_{j,k}^{i}}\right)  \right\vert <\varepsilon,
\]
then $\widetilde{a_{j,k}^{i}}=\tau\left(  \widetilde{a_{j,k}^{i}}\right)
+\left(  \widetilde{a_{j,k}^{i}}-\tau\left(  \widetilde{a_{j,k}^{i}}\right)
\right)  $ and $\left\Vert \left(  \widetilde{a_{j,k}^{i}}-\tau\left(
\widetilde{a_{j,k}^{i}}\right)  \right)  -a_{j,k}^{i}\right\Vert
<2\varepsilon.$ Therefore, no loss of generality, we may assume that
$\widetilde{a_{1,k}^{i}}\in\left(  \mathcal{A}_{\varepsilon}^{1}\right)  ^{0}$
with $\left\Vert a_{1,k}^{i}-\widetilde{a_{1,k}^{i}}\right\Vert <2\varepsilon$
where $k=1,\cdots,n_{i},i=1,\cdots,l.$ And $\widetilde{a_{2,k}^{i}}\in\left(
\mathcal{A}_{\varepsilon}^{2}\right)  ^{0}$ with $\left\Vert a_{2,k}%
^{i}-\widetilde{a_{2,k}^{i}}\right\Vert <2\varepsilon$ where $k=1,\cdots
,n_{i},i=1,\cdots,l.$ Let $\widetilde{b}=\alpha I+\sum_{i=1}^{l}%
\widetilde{a_{1,1}^{i}}\widetilde{a_{2,1}^{i}}\widetilde{a_{1,2}^{i}}%
\cdots\widetilde{a_{1,n_{i}}^{i}}\in\mathcal{A}_{\varepsilon}^{1}%
\ast_{\mathbb{C}}\mathcal{A}_{\varepsilon}^{2}.$ There is an integer $M_{b}>0$
such that
\begin{align*}
&  \left\Vert b-\widetilde{b}\right\Vert \\
&  =\left\Vert \alpha I+\sum_{i=1}^{l}a_{1,1}^{i}a_{2,1}^{i}a_{1,2}^{i}%
a_{2,2}^{i}a_{1,3}^{i}\cdots a_{1,n_{i}}^{i}-\left(  \alpha I+\sum_{i=1}%
^{l}\widetilde{a_{1,1}^{i}}\widetilde{a_{2,1}^{i}}\widetilde{a_{1,2}^{i}%
}\cdots\widetilde{a_{1,n_{i}}^{i}}\right)  \right\Vert \\
&  \leq M_{b}\varepsilon.
\end{align*}
Since $\mathcal{A}_{\varepsilon}^{1}\ast_{\mathbb{C}}\mathcal{A}_{\varepsilon
}^{2}$ is an RFD C*-algebra by Lemma \ref{13.4}, then by Theorem \ref{ttt} we
have that $\mathcal{A}_{1}\ast_{\mathbb{C}}\mathcal{A}_{2}$ is inner QD.
\end{proof}

Since every strong NF algebra is inner, then from \cite{BK1}, we know that
every AF algebra and AH algebra are inner. Hence we have the following two corollaries.

\begin{corollary}
Suppose $\mathcal{A}$ and $\mathcal{B}$ are both AF algebras, then
$\mathcal{A}\ast_{\mathbb{C}}\mathcal{B}$ is an inner QD algebra.
\end{corollary}

\begin{corollary}
Suppose $\mathcal{A}$ and $\mathcal{B}$ are both AH algebras, then
$\mathcal{A}\ast_{\mathbb{C}}\mathcal{B}$ is an inner QD algebra.
\end{corollary}

What will happen when the above amalgamation is over some other C*-algebras
instead of $\mathbb{C}I$? In \cite{LS}, it has been shown that a full
amalgamated free product of two QD algebras may not be MF again, even for a
unital full free product of two full matrix algebras with amalgamation over a
two dimensional C*-algebra which is *-isomorphic to $\mathbb{C\oplus C}$.
Therefore, a unital full amalgamated free product of two unital inner QD
C*-algebras may not be inner again. But we can give the affirmative answers
for some specific cases.

The following result can be found in \cite{P} or \cite{LS}.

\begin{lemma}
\label{32}Suppose that $\mathcal{A},$ $\mathcal{B}$ and $\mathcal{D}$ are
unital C*-algebras. Then
\[
\left(  \mathcal{A\otimes}_{\max}\mathcal{D}\right)  \underset{\mathcal{D}%
}{\mathcal{\ast}}\left(  \mathcal{B\otimes}_{\max}\mathcal{D}\right)
\mathcal{\cong}\left(  \mathcal{A}\underset{\mathbb{C}}{\mathcal{\ast}%
}\mathcal{B}\right)  \otimes_{\max}\mathcal{D}.
\]

\end{lemma}

\begin{lemma}
\label{12f}(\cite{BK2})Let $\mathcal{A}$ be a C*-algebra. Then, for any $k$,
$\mathcal{A}$ is inner QD if and only if $\mathcal{M}_{k}\left(
\mathcal{A}\right)  =\mathcal{A}\otimes\mathcal{M}_{k}\left(  \mathbb{C}%
\right)  $ is inner QD.
\end{lemma}

\begin{proposition}
\label{5}Let $\mathcal{A}$ and $\mathcal{B}$ be unital C*-algebras. If
$\mathcal{D}$ can be embedded as an unital C*-subalgebra of $\mathcal{A}$ and
$\mathcal{B}$ respectively, and $\mathcal{D}$ is *-isomorphic to a full matrix
algebra $\mathcal{M}_{n}\left(  \mathbb{C}\right)  $ for some integer $n,$
then the unital full amalgamated free product $\mathcal{A}\underset
{\mathcal{D}}{\mathcal{\ast}}\mathcal{B}$ is inner QD if $\mathcal{A}$ and
$\mathcal{B}$ are both inner QD.
\end{proposition}

\begin{proof}
Since $\mathcal{D}$ is *-isomorphic to a full matrix algebra$,$ from Lemma
6.6.3 in \cite{KR}, it follows that $\mathcal{A\cong A}^{\prime}%
\otimes\mathcal{D}$ and $\mathcal{B\cong B}^{\prime}\otimes\mathcal{D}$. Then
$\mathcal{A}^{\prime}$ and $\mathcal{B}^{\prime}$ are inner QD by Lemma
\ref{12f}. So the desired conclusion follows from Theorem \ref{12g}, Lemma
\ref{32} and Lemma \ref{12f}.
\end{proof}

Next, we will consider the case when the free products are amalgamated over
some finite-dimensional C*-algebras.

\begin{lemma}
\label{12k}(\cite{BK3}) An arbitrary inductive limit (with injective
connecting maps) of inner quasidiagonal C*-algebras is inner quasidiagonal.
\end{lemma}

\begin{lemma}
(\cite{P}, Theorem 4.2) Assume that we have embeddings of C*-algebras
$\mathcal{C\subseteq A}_{1}\subseteq\mathcal{A}_{2}$ and $\mathcal{C\subseteq
B}_{1}\subseteq\mathcal{B}_{2}$, then the natural morphism $\sigma
:\mathcal{A}_{1}\ast_{\mathcal{C}}\mathcal{B}_{1}\longrightarrow
\mathcal{A}_{2}\ast_{\mathcal{C}}\mathcal{B}_{2}$ is injective.
\end{lemma}

\begin{lemma}
\label{12j}(\cite{P} corollary 4.13) If $\left(  \mathcal{A}_{n}\right)  $ and
$\left(  \mathcal{B}_{n}\right)  $ are increasing sequences of C*-algebras,
all of which contain a common C*-subalgebra $\mathcal{C}$, then there is a
natural isomorphism
\[
\lim_{\longrightarrow}\left(  \mathcal{A}_{n}\ast_{\mathcal{C}}\mathcal{B}%
_{n}\right)  =\lim_{\longrightarrow}\mathcal{A}_{n}\ast_{\mathcal{C}}%
\lim_{\longrightarrow}\mathcal{B}_{n}%
\]
where $\underrightarrow{\lim}$ denotes the ordinary direct limit.
\end{lemma}

The following lemma is a well-known property of AF algebras, we can find it in
\cite{D}

\begin{lemma}
\label{12h}A C*-algebra $\mathcal{A}$ is AF if and only if it is separable and

$\left(  \ast\right)  $ for all $\varepsilon>0$ and $A_{1},\cdots,A_{n}$ in
$\mathcal{A}$, there exists a finite dimensional C*-subalgebra $\mathcal{B}$
of $\mathcal{A}$ such that $dist\left(  A_{i},\mathcal{B}\right)
<\varepsilon$ for $1\leq i\leq n.$

Moreover, if $\mathcal{A}_{1}$ is a finite-dimensional subalgebra of
$\mathcal{A}$, then we may choose $\mathcal{B}$ so that it contains
$\mathcal{A}_{1}.$
\end{lemma}

\begin{lemma}
\label{12l}(\cite{ADRL}, Theorem 4.2) Consider unital inclusions of
C*-algebras $\mathcal{A\supseteq C\subseteq B}$ with $\mathcal{A}$ and
$\mathcal{B}$ finite dimensional. Let $\mathcal{A\ast}_{\mathcal{C}%
}\mathcal{B}$ be the corresponding full amalgamated free product. Then
$\mathcal{A\ast}_{\mathcal{D}}\mathcal{B}$ is RFD if and only there are
faithful tracial states $\tau_{\mathcal{A}}$ on $\mathcal{A}$ and
$\tau_{\mathcal{B}\text{ }}$on $\mathcal{B}$ with
\[
\tau_{\mathcal{A}}(x)=\tau_{\mathcal{B}}(x),\qquad\forall\ x\in\mathcal{C}.
\]

\end{lemma}

\begin{corollary}
Suppose $\mathcal{A}$ and $\mathcal{B}$ are AF algebras and
$\mathcal{A\supseteq C\subseteq B}$ with $\mathcal{C}$ finite-dimensional. If
there are faithful tracial states $\tau_{\mathcal{A}}$ and $\tau_{\mathcal{B}%
}$ on $\mathcal{A}$ and $\mathcal{B}$ respectively, such that
\[
\tau_{\mathcal{A}}(x)=\tau_{\mathcal{B}}(x),\qquad\forall\ x\in\mathcal{C},
\]
then $\mathcal{A}\ast_{\mathcal{C}}\mathcal{B}$ is inner QD..
\end{corollary}

\begin{proof}
Since $\mathcal{C}$ is a finite-dimensional C*-subalgebra, then we can find a
sequence of finite-dimensional C*-subalgebras $\left\{  \mathcal{A}%
_{n}\right\}  _{n=1}^{\mathcal{1}}$ and $\left\{  \mathcal{B}_{n}\right\}
_{n=1}^{\mathcal{1}}$ such that $\mathcal{C\subseteq A}_{1}\subseteq
\mathcal{A}_{2}\subseteq\mathcal{\cdots}$ with $\overline{\cup\mathcal{A}_{n}%
}=\mathcal{A}$ and $\mathcal{C\subseteq B}_{1}\subseteq\mathcal{B}%
_{2}\subseteq\mathcal{\cdots}$ with $\overline{\cup\mathcal{B}_{n}%
}=\mathcal{B}$ by Lemma \ref{12h}. Note that $\mathcal{A}_{n}\ast
_{\mathcal{C}}\mathcal{B}_{n}$ is RFD by Lemma \ref{12l}, then $\mathcal{A\ast
}_{\mathcal{C}}\mathcal{B=}\underrightarrow{\lim}\left(  \mathcal{A}_{n}%
\ast_{\mathcal{C}}\mathcal{B}_{n}\right)  $ is inner by Lemma \ref{12j} and
Lemma \ref{12k}.
\end{proof}

\end{document}